\documentclass[onecolumn]{elsart3p}

\usepackage{graphicx}
\usepackage{color,epsfig}

\usepackage{amssymb}
\usepackage{amsmath}
\usepackage[latin1]{inputenc}
\usepackage{psfrag}

\usepackage[english,francais]{babel}

\newtheorem{theorem}{Theorem}[section]
\newtheorem{lemma}[theorem]{Lemma}
\newtheorem{e-proposition}[theorem]{Proposition}

\newtheorem{e-definition}[theorem]{Definition\rm}

\newtheorem{hyp}[theorem]{Assumption}
\newtheorem{theoreme}{Th\'eor\`eme}[section]

\newtheorem{proposition}[theoreme]{Proposition}

\def\og{\leavevmode\raise.3ex\hbox{$\scriptscriptstyle\langle\!\langle$~}}
\def\fg{\leavevmode\raise.3ex\hbox{~$\!\scriptscriptstyle\,\rangle\!\rangle$}}

\newcommand{\eps}{\varepsilon}
\newcommand{\dps}{\displaystyle}
\newcommand{\RR}{\mathbb R}
\newcommand{\R}{\mathbb R}
\newcommand{\ZZ}{\mathbb Z}
\newcommand{\N}{\mathbb N}
\newcommand{\NN}{\mathbb N}
\newcommand{\E}{\mathbb E}
\newcommand{\cM}{\mathcal{M}}
\newcommand{\cD}{\mathcal{D}}
\newcommand{\cJ}{\mathcal{J}}
\newcommand{\cA}{\mathbb{A}}
\newcommand{\wlim}{\rightharpoonup}

\journal{the Acad\'emie des sciences}
\begin{document}
\centerline{Calcul Scientifique}
\begin{frontmatter}


\selectlanguage{english}
\title{An embedded corrector problem to approximate the homogenized coefficients of an elliptic equation\thanksref{actions}}
\thanks[actions]{We thank Paul Cazeaux for fruitful discussions on questions related to this project. The work of FL is partially supported by ONR under Grant N00014-12-1-0383 and EOARD under grant FA8655-13-1-3061.}


\selectlanguage{english}
\author[authorlabel1,authorlabel4]{Eric Canc\`es},
\ead{cances@cermics.enpc.fr}
\author[authorlabel1,authorlabel4]{Virginie Ehrlacher},
\ead{ehrlachv@cermics.enpc.fr}
\author[authorlabel2,authorlabel4]{Fr\'ed\'eric Legoll},
\ead{legoll@lami.enpc.fr}
\author[authorlabel3]{Benjamin Stamm}
\ead{stamm@ann.jussieu.fr}

\address[authorlabel1]{CERMICS, \'Ecole des Ponts ParisTech, 77455 Marne-La-Vall\'ee Cedex 2, France}
\address[authorlabel2]{Laboratoire Navier, \'Ecole des Ponts ParisTech, 77455 Marne-La-Vall\'ee Cedex 2, France}
\address[authorlabel3]{Sorbonne Universit\'es, UPMC Univ. Paris 06 and CNRS, UMR 7598, Laboratoire J.-L. Lions, 75005 Paris, France}
\address[authorlabel4]{INRIA Rocquencourt, MATHERIALS research-team, Domaine de Voluceau,  B.P. 105, 78153 Le Chesnay Cedex, France}


\medskip
\begin{center}
{\small Received *****; accepted after revision +++++\\
Presented by *****}
\end{center}

\begin{abstract}
\selectlanguage{english}
We consider a diffusion equation with highly oscillatory coefficients that admits a homogenized limit. As an alternative to standard corrector problems, we introduce here an embedded corrector problem, written as a diffusion equation in the whole space in which the diffusion matrix is uniform outside some ball of radius $R$. Using that problem, we next introduce three approximations of the homogenized coefficients. These approximations, which are variants of the standard approximations obtained using truncated (supercell) corrector problems, are shown to converge when $R \to \infty$. We also discuss efficient numerical methods to solve the embedded corrector problem.

\vskip 0.5\baselineskip

\selectlanguage{francais}
\noindent{\bf R\'esum\'e} \vskip 0.5\baselineskip \noindent
{\bf Un probl\`eme de correcteur incorpor\'e pour approcher les coefficients homog\'en\'eis\'es d'une \'equation elliptique.}

Nous consid\'erons une \'equation de diffusion \`a coefficients hautement oscillants qui admet une limite homog\'en\'eis\'ee, et nous introduisons une variante du probl\`eme du correcteur standard, que nous appelons probl\`eme du correcteur incorpor\'e. Celui-ci s'\'ecrit comme une \'equation de diffusion pos\'ee dans tout l'espace, dans laquelle la matrice de diffusion est uniforme \`a l'ext\'erieur d'une boule de rayon $R$. Nous introduisons ensuite trois approximations des coefficients homog\'en\'eis\'es, calcul\'ees \`a partir de la solution de ce probl\`eme. Ces approximations, qui sont des variantes des approximations standard bas\'ees sur le probl\`eme du correcteur tronqu\'e (m\'ethode de supercellule), convergent lorsque $R \to \infty$. Nous mentionnons \'egalement des m\'ethodes de r\'esolution num\'erique efficaces du probl\`eme du correcteur incorpor\'e.

\end{abstract}
\end{frontmatter}

\selectlanguage{english}
\section{Introduction}

We consider the standard elliptic, highly oscillatory problem
\begin{equation}
  \label{eq:div-eps}
  -\hbox{\rm div}\,\left[ A\left(\cdot/\varepsilon\right) \, \nabla u_\varepsilon \right]=f \ \text{on $\Omega$}, \quad u_\eps = 0 \ \text{on $\partial \Omega$}, 
\end{equation}
where $\Omega$ is a smooth bounded domain of $\RR^d$ and $f\in L^2(\Omega)$. The coefficient~$A$ is a matrix-valued field, and $\varepsilon$ is a small characteristic length-scale. Throughout this Note, we assume that $A$ is symmetric and elliptic, in the sense that there exists $0 < \alpha \leq \beta < \infty$ such that $A(x) \in \cM_{\alpha, \beta}$ for any $x \in \R^d$, where
$$
\cM_{\alpha, \beta} := \left\{ A\in \R^{d\times d}, \; A^T = A \ \text{and, for any $\xi \in \RR^d$}, \ \alpha |\xi|^2 \leq \xi^T A \xi \leq \beta |\xi|^2 \right\}.
$$
It is well-known (see e.g.~\cite{blp,cd,jikov}) that, under this assumption, problem~\eqref{eq:div-eps} admits a homogenized limit, i.e. that the sequence $A(\cdot/\eps)$ $G$-converges, up to the extraction of a subsequence, to some homogenized matrix-valued field $A^\star \in L^\infty(\Omega,\cM_{\alpha, \beta})$ when $\eps \to 0$ (the notion of $G$-convergence is recalled in Definition~\ref{def:G-conv} below). 

Our setting includes in particular the periodic case, where $A(x)=A_{\rm per}(x)$ for a fixed $\ZZ^d$-periodic function $A_{\rm per}$, and the random stationary case (see~\cite{papa}), where
\begin{equation}
\label{eq:cas_random}
A(x)=A_{\rm sta}(x,\omega) \ \text{for some realization $\omega$ of a random stationary function $A_{\rm sta}$.}
\end{equation}
In these two cases, the whole sequence $A(\cdot/\eps)$ $G$-converges (for almost all $\omega$ in the case~\eqref{eq:cas_random}).

\medskip

Computing the homogenized coefficient $A^\star$ is in general a challenging task, even in the cases when a closed form formula for $A^\star$ is available. Consider for instance the random stationary case~\eqref{eq:cas_random} in a discrete stationary setting~\cite{singapour,jmpa} when 
$$
\forall k \in \ZZ^d, \quad A_{\rm sta}(x,\tau_k \omega) = A_{\rm sta}(x+k,\omega) \quad \text{a.e. in $x$, a.s. in $\omega$},
$$ 
where $(\tau_k)_{k\in \ZZ^d}$ is an ergodic group action on the probability space. In that setting, $A^\star$ is a constant deterministic matrix, given by
\begin{equation}
\label{eq:defAstar}
\forall p \in \R^d, \quad A^\star p = \E \left[\frac{1}{|Q|} \int_Q A(x, \cdot) \left( p + \nabla w_p(x, \cdot) \right) \,dx \right], \quad Q=(0,1)^d,
\end{equation}
where $w_p$ is the unique solution (up to an additive constant) to the so-called corrector problem
\begin{equation}
\label{eq:correctorrandom}
\left\{
\begin{array}{l}
-\hbox{\rm div}\left[ A(\cdot, \omega)(p + \nabla w_p(\cdot, \omega))\right] = 0 \hbox{ almost surely in $\cD'(\R^d)$},\\
\dps \nabla w_p \hbox{ is stationary}, \quad \E\left[ \int_Q \nabla w_p(x, \cdot)\,dx \right] = 0.
\end{array}
\right.
\end{equation}
The major difficulty to compute $A^\star$ is the fact that the corrector problem~\eqref{eq:correctorrandom} is set over the whole space $\R^d$ and cannot be reduced to a problem posed over a bounded domain (in contrast to e.g. periodic homogenization). This is the reason why approximation strategies are required, yielding practical approximations of $A^\star$. A popular approach, introduced in~\cite{bourgeat}, is to approximate $A^\star$ by $A^\star_N(\omega)$, which, in turn, is defined by
\begin{equation}
\label{eq:defAstar_N}
\forall p \in \R^d, \quad A^\star_N(\omega) p:= \frac{1}{|Q_N|} \int_{Q_N} A(x, \omega) \left( p + \nabla w^N_p(x, \omega) \right) \,dx, \quad Q_N=(-N,N)^d,
\end{equation}
where $w_p^N$ is the unique solution (up to an additive constant) to the truncated corrector problem
\begin{equation}
\label{eq:correctorrandom_N}
-\hbox{\rm div}\left[ A(\cdot, \omega)(p + \nabla w^N_p(\cdot, \omega))\right] = 0 \hbox{ almost surely in $\cD'(\R^d)$},
\quad
\hbox{$w^N_p(\cdot,\omega)$ is $Q_N$-periodic}.
\end{equation}
As shown in~\cite{bourgeat}, $A^*_N(\omega)$ almost surely converges to $A^\star$ when $N \to \infty$.

\medskip

The aim of this Note is to introduce variants of~\eqref{eq:defAstar_N}--\eqref{eq:correctorrandom_N} that allow to compute accurate approximations of the homogenized coefficient $A^\star$, and that, in some cases, are amenable to efficient numerical implementations through the use of boundary integral formulations. We refer to~\cite{cras_kun_li} for other characterizations of the homogenized matrix, which can also be turned into numerical strategies alternative to~\eqref{eq:defAstar_N}--\eqref{eq:correctorrandom_N} to approximate $A^\star$ in the random stationary setting. See also~\cite{cottereau,gloria-esaim} for other numerical strategies to approximate~\eqref{eq:defAstar}.

In Section~\ref{sec:embedded}, we describe our approach and explain in what sense it is amenable to an efficient implementation. Based on that approach, alternative approximations of $A^\star$ are built in Section~\ref{sec:new_defs}, where we also collect convergence results. The results presented in this Note will be complemented and extended in~\cite{futur}.

\section{Embedded corrector problem}
\label{sec:embedded}

In this section, we introduce an embedded corrector problem (see~\eqref{eq:pbbase} below), which is key in our approach.

\medskip

For any $R>0$, we denote by $B_R$ the open ball of $\R^d$ of radius $R$ centered at the origin, and $B := B_1$. Let $\Gamma_R:= \partial B_R$ and $n_R(x)$ be the normal unitary vector of $\Gamma_R$ at the point $x \in \Gamma_R$ pointing outwards $B_R$. We introduce the vector spaces
$$
V:=\left\{v\in L^2_{\rm loc}(\R^d), \; \nabla v \in \left(L^2(\R^d)\right)^d\right\} \quad \hbox{ and } \quad V_0:= \left\{ v \in V, \; \int_B v = 0\right\}.
$$
The space $V_0$, endowed with the scalar product $\langle \cdot, \cdot \rangle$ defined by
$$
\forall v,w\in V_0, \quad \langle v,w\rangle:= \int_{\R^d} \nabla v \cdot \nabla w,
$$
is a Hilbert space. 

\medskip

\begin{figure}[htbp]
\psfrag{A}{$A_{\rm ext}$}
\psfrag{B}{$A_{\rm int}^i$}
\psfrag{C}{$\partial B_R$}
\psfrag{D}{$A$}
\psfrag{E}{$\partial B_R$}
\centerline{
\includegraphics[width=6cm]{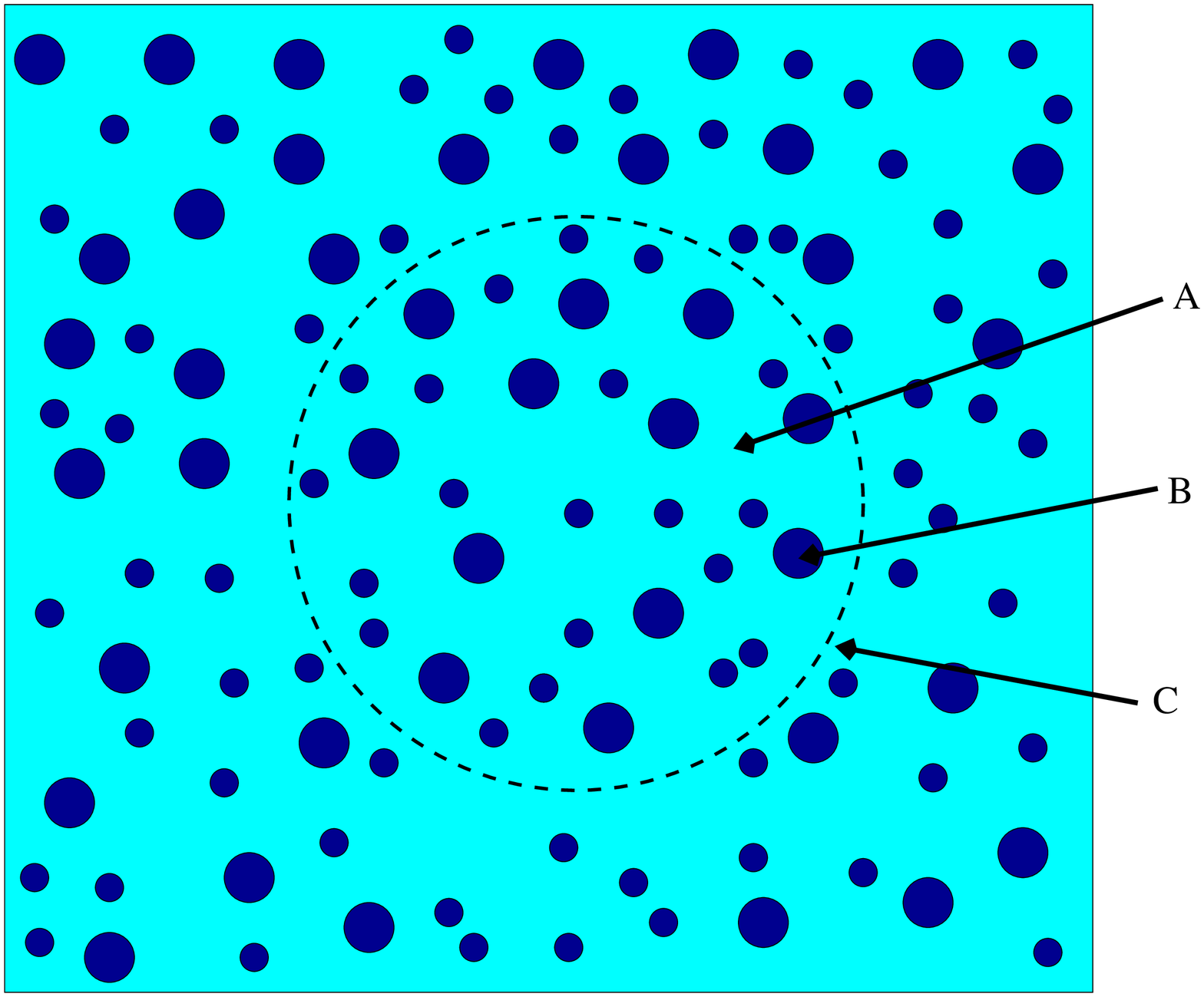}
\qquad \qquad
\includegraphics[width=6cm]{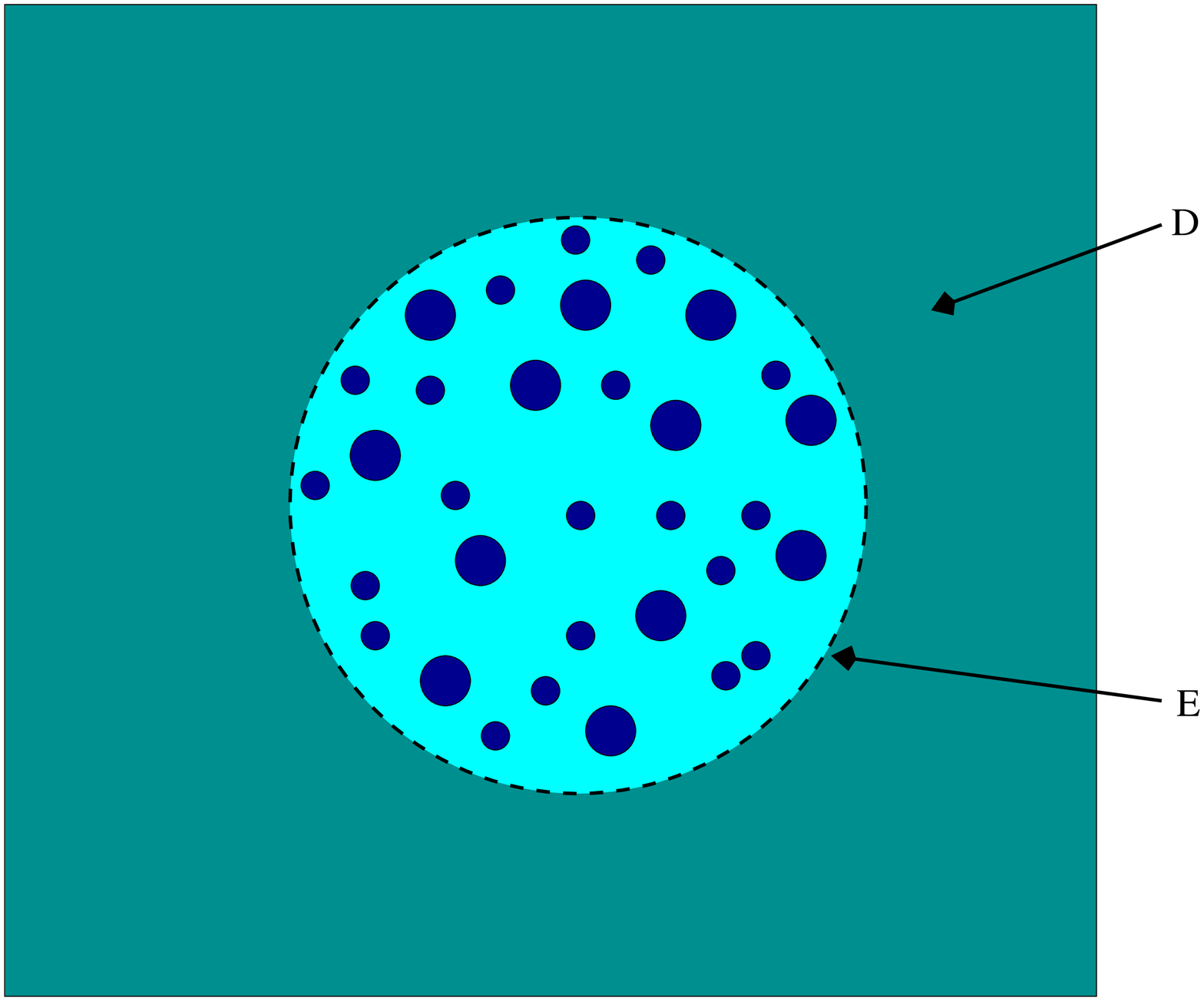}
}
\caption{Left: field $\cA(x)$. Right: field $\cA_{R,A}(x)$: beyond the sphere of radius $R$, the field $\cA(x)$ is replaced by a uniform coefficient~$A$.\label{fig:boules}}
\end{figure}

For any matrix-valued field $\cA \in L^\infty(\R^d,\cM_{\alpha,\beta})$, any $R>0$, any constant matrix $A\in \cM_{\alpha,\beta}$, and any vector $p\in \R^d$, we denote by $w^{R,\cA,A}_p$ the unique solution in $V_0$ to 
\begin{equation}
\label{eq:pbbase}
\fbox{$
-\hbox{\rm div} \Big( \cA_{R,A} (p + \nabla w^{R,\cA, A}_p)\Big) = 0 \hbox{ in } \cD'(\R^d),
$}
\end{equation}
where (see Figure~\ref{fig:boules})
$$
\cA_{R,A}(x) := \left| 
\begin{array}{l}
\cA(x) \mbox{ if } x\in B_R,\\
A \mbox{ if } x\in \R^d \setminus B_R.
\end{array}\right .
$$
In~\eqref{eq:pbbase}, we keep the original coefficient $\cA$ in the ball $B_R$, and replace it outside $B_R$ by a uniform coefficient~$A$.

\medskip

Assume that the matrix-valued field $\cA \in L^\infty(\R^d,\cM_{\alpha, \beta})$ satisfies the following:

\medskip

\begin{hyp}
\label{hyp:G-conv}
The rescaled matrix-valued fields $\cA^R$, defined by $\cA^R(x)=\cA(Rx)$, form a family $(\cA^R)_{R>0}$ that $G$-converges to a constant matrix $A^\star \in \cM_{\alpha, \beta}$ on $B$ as $R$ tends to infinity.
\end{hyp}

\medskip

Under this assumption, the motivation for considering problems of the form~\eqref{eq:pbbase} is twofold. First, we show in Section~\ref{sec:new_defs} below that the solution $w^{R,\cA, A}_p$ to~\eqref{eq:pbbase} can be used to define consistent approximations of $A^\star$. Second, in some cases, problem~\eqref{eq:pbbase} can be efficiently solved, using a numerical approach directly inspired from that proposed in~\cite{Stamm1,Stamm2}. This is for example the case when, in $B_R$,
\begin{equation}
\label{eq:spherical}
\cA(x)=\left|
\begin{array}{l}
A_{\rm int}^i \mbox{ if } x \in B_R \cap B(x_i,r_i), \ 1 \leq i \leq I,
\\
A_{\rm ext} \mbox{ if } x \in B_R \setminus \bigcup_{i=1}^I B(x_i, r_i),
\end{array}
\right.
\end{equation}
for some $I\in\N^\star$, $A_{\rm int}^i, A_{\rm ext} \in \cM_{\alpha, \beta}$ for any $1 \leq i \leq I$, $(x_i)_{1\leq i \leq I}\subset B_R$ and $(r_i)_{1\leq i \leq I}$ some set of positive real numbers such that $\bigcup_{i=1}^I B(x_i,r_i) \subset B_R$ and $B(x_i,r_i) \cap B(x_j,r_j) = \emptyset$ for all $1\leq i \neq j \leq I$. We have denoted by $B(x_i,r_i) \subset \R^d$ the ball of radius $r_i$ centered at $x_i$. We refer to~\cite{futur} for other cases.

The expression~\eqref{eq:spherical} corresponds to the case of (possibly stochastic) heterogeneous materials composed of spherical inclusions. The properties of the inclusions (i.e. the coefficients $A_{\rm int}^i$), their centers $x_i$ and their radii $r_i$ may be random, as long as $\cA$ is stationary (see Figure~\ref{fig:boules}).

In the case~\eqref{eq:spherical}, problem~\eqref{eq:pbbase} can be efficiently solved using a boundary integral method (see~\cite{futur}). Since $\cA_{R,A}$ is uniform in each $B(x_i,r_i)$, in $\dps B_R \setminus \cup_i B(x_i,r_i)$ and in $\R^d \setminus B_R$, problem~\eqref{eq:pbbase} can indeed be recast as an integral equation on the spheres $\partial B(x_i,r_i)$ and $\Gamma_R$. In the case of random homogenization, the practical consequence is that, for the same number of degrees of freedom, we can afford to work on domains $B_R$ that are much larger than the truncated domains $Q_N$ in~\eqref{eq:defAstar_N}--\eqref{eq:correctorrandom_N}. We thus expect to obtain better approximations of $A^\star$.

\section{New definitions of approximate homogenized matrices} 
\label{sec:new_defs}

Assume that the matrix field $\cA$ satisfies Assumption~\ref{hyp:G-conv}. We wish to use solutions to~\eqref{eq:pbbase} to construct a family $(A^{\star,R})_{R>0}$ which converges to the homogenized matrix $A^\star$ as $R$ tends to infinity.

In the subsequent Sections~\ref{sec:def1}, \ref{sec:def2} and \ref{sec:def3}, we respectively present three possible choices leading to converging approximations, namely~\eqref{eq:optimisation}, \eqref{eq:def2} and~\eqref{eq:def3}. We refer to~\cite{futur} for the proof of the results stated below. To introduce these choices, we note that the solution $w_p^{R,\cA, A}$ to~\eqref{eq:pbbase} is equivalently the unique solution to the optimization problem
$$
w_p^{R,\cA, A} = \mathop{\mbox{argmin}}_{v\in V_0} J_p^{R,\cA, A}(v), 
$$
where
\begin{equation}
\label{eq:optim2}
J^{R,\cA, A}_p(v) := \frac{1}{2|B_R|} \int_{B_R} (p + \nabla v)^T \cA (p + \nabla v) + \frac{1}{2|B_R|}\int_{\R^d\setminus B_R} (\nabla v)^T A \nabla v - \frac{1}{|B_R|} \int_{\Gamma_R} (Ap\cdot n_R) v.
\end{equation}
We set
$$
\cJ_p^{R,\cA}(A) := J_p^{R,\cA, A}(w_p^{R,\cA, A}) = \min_{v\in V_0} J_p^{R,\cA, A}(v). 
$$
The linearity of the mapping $\R^d \ni p \mapsto w_p^{R,\cA, A} \in V_0$ yields that, for any $A\in \cM_{\alpha, \beta}$, there exists a unique symmetric matrix $G^{R,\cA}(A) \in \RR^{d \times d}$ such that
\begin{equation}
\label{eq:defGA}
\fbox{$\dps
\forall p\in \R^d, \quad \frac{1}{2}p^T G^{R,\cA}(A) p = \cJ^{R,\cA}_p(A).
$}
\end{equation}
Note that $\dps \frac{1}{2}\mbox{Tr}\left(G^{R,\cA}(A)\right) = \sum_{i=1}^d \cJ^{R,\cA}_{e_i}(A)$, where $(e_i)_{1 \leq i \leq d}$ is the canonical basis of $\RR^d$. The following expression of $\cJ_p^{R,\cA}(A)$ is useful:
\begin{multline}
\label{eq:expJ3}
\cJ_p^{R,\cA}(A) = \frac{1}{2|B_R|}\int_{B_R} p^T \cA p - \frac{1}{2|B_R|} \int_{B_R} \left( \nabla w^{R,\cA, A}_p\right)^T \cA \nabla w^{R,\cA, A}_p 
\\
- \frac{1}{2|B_R|} \int_{\R^d\setminus B_R} \left( \nabla w^{R,\cA, A}_p\right)^T A \nabla w^{R,\cA, A}_p. 
\end{multline}
Before describing our three approaches, we recall the following classical definition (see~\cite{MuratTartar}):

\medskip

\begin{e-definition}[$G$-convergence]
\label{def:G-conv}
Let $D$ be an open bounded smooth subdomain of $\R^d$. A family of matrix-valued functions $\left( \cA^R \right)_{R>0} \subset L^\infty(D, \cM_{\alpha, \beta})$ is said to $G$-converge in $D$ to a matrix-valued function $\cA^\star\in L^\infty(D, \cM_{\alpha, \beta})$ if, for all $f\in H^{-1}(D)$, the family $(u^R)_{R>0}$ of solutions to 
$$
-\mbox{\rm div}\left( \cA^R \nabla u^R \right) = f \mbox{ in } \cD'(D), \quad u^R\in H^1_0(D),
$$
satisfies
$$
u^R \mathop{\wlim}_{R\to +\infty} u^\star \mbox{ weakly in }H^1_0(D),
\qquad
\cA^R \nabla u^R \mathop{\wlim}_{R\to +\infty} \cA^\star \nabla u^\star \mbox{ weakly in }L^2(D),
$$
where $u^\star$ is the unique solution to the homogenized equation
$$
-\mbox{\rm div}\left( \cA^\star \nabla u^\star \right) = f \mbox{ in } \in \cD'(D), \quad u^\star \in H^1_0(D).
$$
\end{e-definition}

\subsection{First alternative definition: minimizing the scattering energy}
\label{sec:def1}

To gain some intuition, we first recast~\eqref{eq:pbbase} as
$$
-\mbox{\rm div}\left[ \Big(A + \chi_{B_R} (\cA -A) \Big) \, \Big(p + \nabla w_p^{R,\cA, A} \Big)\right] = 0 \quad \mbox{ in } \quad \cD'(\R^d),
$$
where $\chi_{B_R}$ is the characteristic function of $B_R$. Thus, in this problem, the quantity $\cA-A$ can be seen as a local perturbation to the homogeneous exterior medium characterized by the diffusion coefficient $A$. In turn, $w_p^{R,\cA, A}$ can be seen as the perturbation of an incident plane wave with wavevector $p$ induced by the defect located in $B_R$. This is somehow reminiscent of the classical Eshelby problem~\cite{Eshelby}. A first idea is to choose a constant exterior matrix such that the scattering energy of the perturbation of the wave is as small as possible. We have the following result (recall that $G^{R,\cA}$ is defined by~\eqref{eq:defGA}):

\medskip

\begin{lemma}
\label{lem:lemconcavity}
For all $R>0$ and $\cA\in L^\infty(\R^d,\cM_{\alpha, \beta})$, the function $\dps \cM_{\alpha, \beta}\ni A \mapsto \mbox{Tr}\left(G^{R,\cA}(A)\right)$ is concave. 
\end{lemma}

\medskip

It follows that, for any $R>0$, there exists (at least) one matrix $A^R_1 \in \cM_{\alpha, \beta}$ such that
\begin{equation}
\label{eq:optimisation}
\fbox{$\dps
A^R_1 = \mathop{\mbox{argmax}}_{A\in \cM_{\alpha, \beta}} \ \ \mbox{Tr}\left(G^{R,\cA}(A)\right).
$}
\end{equation}
This matrix $A^R_1$ can be seen as a matrix which minimizes the scattering energy induced by the defect $\cA - A$ of incident plane waves in an infinite medium. Indeed, using~\eqref{eq:expJ3}, we have that
$$
A^R_1 = \mathop{\mbox{argmin}}_{A\in \cM_{\alpha, \beta}} \sum_{i=1}^d \left( \int_{B_R} \left( \nabla w_{e_i}^{R,\cA, A}\right)^T \cA \nabla w_{e_i}^{R,\cA, A} + \int_{\R^d \setminus B_R} \left( \nabla w_{e_i}^{R,\cA, A}\right)^T A \nabla w_{e_i}^{R,\cA, A} \right),
$$
and the matrix $A^R_1$ can thus be seen as a diffusion matrix $A$ of the exterior medium such that the sum of the energies of the scattering waves with incident wavevectors $e_i$ induced by the defect is minimum.

\medskip

As shown in Proposition~\ref{prop:prop1} below, the approximation $A^R_1$ converges to $A^\star$ when $R \to \infty$.

\subsection{Second alternative definition: an equivalent internal homogeneous material}
\label{sec:def2}

We now introduce a second alternative definition of an approximate homogenized matrix:
\begin{equation}
\label{eq:def2}
\fbox{$
A_2^R = G^{R,\cA}(A^R_1),
$}
\end{equation}
where $A^R_1$ is defined by~\eqref{eq:optimisation}. In view of~\eqref{eq:defGA}, the above relation can also be written as
$$
\forall p\in\R^d, \quad \frac{1}{2} p^T A_2^R p = \cJ_{p}^{R,\cA}(A^R_1).
$$
Using~\eqref{eq:optim2}, the above definition can formally be recast as
\begin{multline}
\label{eq:formal}
\int_{B_R} \left(p + \nabla w_p^{R,\cA, A^R_1}\right)^T \cA \left(p + \nabla w_p^{R,\cA, A^R_1}\right) + \int_{\R^d\setminus B_R} \left(p + \nabla w_p^{R,\cA, A^R_1}\right)^T A^R_1 \left(p + \nabla w_p^{R,\cA, A^R_1}\right) 
\\
= \int_{B_R} p^T A^R_2 p + \int_{\R^d\setminus B_R} p^T A^R_1 p.
\end{multline}
The above relation is formal in the sense that both sides of the equation are infinite, but it nevertheless has an interesting physical interpretation. The above left-hand side corresponds to the energy of the heterogeneous material, modelled by $\cA$ in $B_R$ and $A^R_1$ outside of $B_R$, and where the field $p + \nabla w_p^{R,\cA, A^R_1}$ is solution to the equilibrium equation~\eqref{eq:pbbase}. Since $\nabla w_p^{R,\cA, A^R_1}$ is in $L^2(\R^d)$, its average is thought to vanish, and hence the average field is $p$. The above right-hand side corresponds to the energy of a material, modelled by $A^R_2$ in $B_R$ and $A^R_1$ outside of $B_R$, in which the field is uniform and equal to $p$. The formal equation~\eqref{eq:formal} thus ``defines'' $A^R_2$ by an equality in terms of energies. 

\medskip

The following convergence result can be established:
\begin{e-proposition}
\label{prop:prop1}
Assume that the matrix field $\cA$ satisfies Assumption~\ref{hyp:G-conv}. Then, the two families of matrices $\left( A^R_1 \right)_{R>0}$ and $\left( A^R_2 \right)_{R>0}$, respectively defined by~\eqref{eq:optimisation} and~\eqref{eq:def2}, satisfy
$$
A^R_1 \mathop{\longrightarrow}_{R\to +\infty} A^\star 
\quad \mbox{ and } \quad
A^R_2 \mathop{\longrightarrow}_{R\to +\infty} A^\star. 
$$
\end{e-proposition}

\subsection{Third alternative definition: a self-consistent equation}
\label{sec:def3}

We eventually introduce a third alternative definition, inspired by the approximation of $A^\star$ introduced in~\cite{Christensen}. Assume that, for any $R>0$, there exists a matrix $A^R_3 \in \cM_{\alpha, \beta}$ such that  
\begin{equation}
\label{eq:def3}
\fbox{$
A^R_3  = G^{R,\cA}(A^R_3).
$}
\end{equation}
Such a matrix formally satisfies the self-consistent equation
\begin{multline*}
\sum_{i=1}^d \int_{B_R} \left[\left(e_i + \nabla w_{e_i}^{R,\cA, A^R_3}\right)^T \cA \left(e_i + \nabla w_{e_i}^{R,\cA, A^R_3}\right) - e_i^T A^R_3 e_i \right] 
\\
+ \int_{\R^d\setminus B_R}\left[ \left(e_i + \nabla w_{e_i}^{R,\cA, A^R_3}\right)^T A^R_3 \left(e_i + \nabla w_{e_i}^{R,\cA, A^R_3}\right) - e_i^T A^R_3 e_i \right] = 0. 
\end{multline*}
This third definition also yields a converging approximation of $A^\star$:

\medskip

\begin{proposition}
\label{prop:prop2}
Assume that the matrix field $\cA$ satisfies Assumption~\ref{hyp:G-conv}, and that there exists a sequence $\left( A^{R_k}_3 \right)_{k \in \NN} \in \left( \cM_{\alpha, \beta} \right)^\NN$ satisfying
$$
\forall k \in \NN, \quad A^{R_k}_3  = G^{R_k,\cA}\left(A^{R_k}_3\right)
$$
for some increasing sequence $\left( R_k \right)_{k \in \NN}$ of positive numbers converging to $+\infty$. Then, 
$$
A^{R_k}_3 \mathop{\longrightarrow}_{k\to +\infty} A^\star. 
$$
\end{proposition}

Note that we do not assume in this Proposition that the fixed point equation~\eqref{eq:def3} has a solution for {\em all} radii $R$. Proving the existence of a matrix $A^R_3$ satisfying~\eqref{eq:def3} in the general case is a delicate question. We however already have the following partial result, which addresses the isotropic case.

\medskip

\begin{proposition}
\label{prop:prop3}
Let $d \geq 2$. Let $\cA \in L^\infty(\R^d, \cM_{\alpha, \beta})$ be a matrix-valued field satisfying Assumption~\ref{hyp:G-conv}. Assume also that the homogenized matrix satisfies $A^\star = a^\star I_d$, where $I_d$ is the identity matrix of $\RR^{d \times d}$.

Then $a^\star \in [\alpha,\beta]$ and, for any $R>0$, there exists $a^R_3 \in [\alpha,\beta]$ such that 
\begin{equation}
\label{eq:spher}
a^R_3 = \frac{1}{d}\mbox{Tr}\left( G^{R,\cA}(a^R_3 I_d)\right).
\end{equation}
In addition, 
$$
a^R_3 \mathop{\longrightarrow}_{R\to +\infty} a^\star. 
$$
\end{proposition}
Note that~\eqref{eq:spher} is weaker than~\eqref{eq:def3}, which would read in this case $a^R_3 I_d = G^{R,\cA}(a^R_3 I_d)$. However, this weaker condition is sufficient to prove that $a^R_3$ is a converging approximation of $a^\star$.

\medskip

We conclude with the following two remarks. First, in the one-dimensional case, it is possible to obtain explicit expressions for $A^R_1$, $A^R_2$ and $A^R_3$ (which are uniquely defined by~\eqref{eq:optimisation}, \eqref{eq:def2} and~\eqref{eq:spher}, respectively) and see that they converge to $A^\star$ when $R \to \infty$. Second, in the case when $\cA$ is actually equal to a constant matrix $A$ in $B_R$, then we have $A^R_1 = A^R_2 = A$, while the unique solution to~\eqref{eq:def3} is $A^R_3 = A$.




\end{document}